\documentclass[
               oneeqnum,onefignum
               ]{siamltex}
\usepackage[english]{babel}
\usepackage{color}
\usepackage{hyperref}

\usepackage{amssymb}

\date{}

\newtheorem{remark}{Remark}

\renewenvironment{remark}[1][\hspace{-1.0ex}]%
 {\par\addvspace{0mm}\indent\refstepcounter{remark} {\it Remark.} }%
 {\par\addvspace{0mm}\rm}

\def\laaa{[}
\def\raaa{]}
\def\lbbb{[}
\def\rbbb{]}

\title{$n$-Ary Quasigroups of Order $4$}
\author{\lowercase{\href{http://arxiv.org/find/grp_math/1/au:+Krotov_Denis/0/1/0/all/0/1}{\uppercase{Denis~S.~Krotov}}}%
\thanks{Sobolev Institute of Mathematics,
prosp. Akademika Koptyuga, 4, Novosibirsk, 630090, Russia
({\tt krotov@math.nsc.ru}).
}
\and \lowercase{\href{http://arxiv.org/find/grp_math/1/au:+Potapov_Vladimir/0/1/0/all/0/1}{\uppercase{Vladimir~N.~Potapov}}}%
\thanks{Sobolev Institute of Mathematics,
prosp. Akademika Koptyuga, 4, Novosibirsk, 630090, Russia
({\tt vpotapov@math.nsc.ru}).
}
}

\begin{document}
\maketitle
\begin{abstract}
We characterize the set of all $n$-ary qua\-si\-groups of order $4$:
every $n$-ary qua\-si\-group of order $4$
is permutably reducible or semilinear.
Permutable reducibility means that an $n$-ary qua\-si\-group
can be represented as a composition of $k$-ary and $(n-k+1)$-ary
quasigroups for some $k$ from $2$ to $n-1$,
where the order of arguments in the representation
can differ from the original order.
The set of semilinear $n$-ary quasigroups
has a characterization in terms of Boolean functions.
\end{abstract}

\section{Introduction}

An algebraic system consisting of a finite set $\Sigma$ of cardinality $|\Sigma|=q$
and an $n$-ary operation $f: \Sigma^n\rightarrow \Sigma$
uniquely invertible in each place is called an $n$-ary quasigroup of order $q$.
The function $f$ can also
be referred to as an $n$-ary quasigroup of order $q$ or, for short,
an $n$-qua\-si\-group.
The value array of an $n$-qua\-si\-group of order $q$
is known as a latin $n$-cube of order $q$
(if $n=2$, a latin square).
Furthermore, there is a one-to-one
correspondence between the $n$-quasigroups and the distance $2$ MDS codes in $\Sigma^{n+1}$.

It is known that for every $n$ there exist
exactly two equivalent $n$-qua\-si\-groups of order $2$ and
$3\cdot 2^n$ $n$-qua\-si\-groups of order $3$,
which constitute one isotopy class
(see, e.g., \cite{LaywineMullen}).
So, $4$ is the first order for which a rich class of $n$-qua\-si\-groups exists.
On the other hand, this order is of special interest for different areas
of mathematics close to information theory.
For example,
\begin{remunerate}
\item[$\bullet$] the class of $1$-perfect codes in $\{0,1\}^n$ of rank at most $n-\log_2 (n+1)+2$
(the minimum rank is $n-\log_2 (n+1)$ for $1$-perfect codes)
is characterized in terms of $n$-qua\-si\-groups of order $4$, see \cite{AvgHedSol:class}
(so, our work completes this characterization);
\item[$\bullet$] order $4$ is the first order that is applicable for use in quasigroup stream ciphers;
\item[$\bullet$] from $n$-qua\-si\-groups of order $4$, $n$-qua\-si\-groups of other orders can be constructed,
giving examples of $n$-qua\-si\-groups with nontrivial properties
(see, e.\,g., \linebreak[5] \cite{Kro:n-2}).
\end{remunerate}

In this paper, we show that every $n$-qua\-si\-group of order $4$
is permutably reducible or semilinear.
Permutable reducibility means that the $n$-qua\-si\-group can be represented as a repetition-free
composition of qua\-si\-groups of smaller arities were
the ordering of the arguments
in the representation
can differ from the original (see Definition~\ref{def:pr}).
Semilinearity (Definition~\ref{def:SSL}) means that the $n$-qua\-si\-group can be obtained as a direct product
of two $n$-qua\-si\-groups of order $2$ modified by a Boolean function $\{0,1\}^n\to\{0,1\}$
(sometimes this construction is referred to as the wreath product construction, but we should
remember that this does not agree with the concept of wreath product of groups).

In \S\ref{s:1} we introduce main concepts and notations.
In \S\ref{s:2} we formulate the result (Theorem~\ref{th:main})
and divide the proof into four subcases, Lemmas~\ref{pro:3b}-\ref{pro:1b}.
Lemmas~\ref{pro:4b} and~\ref{pro:1b} are proved in \S\ref{sect:4b} and \S\ref{sect:1b},
while Lemmas~\ref{pro:3b} and~\ref{pro:2b} follow from previous papers.

\section{Main definitions}\label{s:1}\mbox{}

\begin{definition}
{
 An $n$-ary operation $Q:\Sigma^n\to \Sigma$,
where $\Sigma$ is a nonempty set,
is called
an \emph{$n$-ary qua\-si\-group} or \emph{$n$-qua\-si\-group}
$($\emph{of order $|\Sigma|$}$)$ if 
in the equality $z_{0}=Q(z_1, \ldots , z_n)$ knowledge of any $n$ elements
of $z_0$, $z_1$, \ldots , $z_n$ uniquely specifies the remaining one {\rm \cite{Belousov}}.
}
\end{definition}

The definition is symmetric with respect to the variables
$z_0$, $z_1$, \ldots , $z_n$, and sometimes it is convenient to use
a symmetric form for the relation $z_{0}=Q(z_1, \ldots , z_n)$.
For this reason, we will denote by $Q\langle z_0, z_1, \ldots , z_n\rangle$
the corresponding predicate, i.\,e., the characteristic function of this relation.
(In coding theory, the set corresponding to this predicate is known as a distance $2$ MDS code.)

Given $\bar y = (y_1,\dots,y_n)$, we denote
$$\bar y^{\laaa i\raaa}\lbbb x\rbbb \triangleq (y_1,\ldots,y_{i-1},x,y_{i+1},\ldots,y_n);$$
similarly, we define
$\bar y^{\laaa i_1,i_2,\dots,i_k\raaa}\lbbb x_{i_1},x_{i_2},\ldots,x_{i_k}\rbbb$.

\begin{definition}
{
 If we assign some fixed values to $l\in \{1,\ldots,n\}$ variables
in the predicate $Q\langle z_0, \ldots , z_n\rangle$
then the $(n-l+1)$-ary predicate obtained corresponds to an $(n-l)$-qua\-si\-group.
Such a qua\-si\-group is called a \emph{retract} or \emph{$(n-l)$-retract} of $Q$.
If $z_0$ is not fixed, the retract is \emph{principal}.
}
\end{definition}

\begin{definition}
{
 By an {\em isotopy} we shall mean
a collection of $n+1$ permutations $\tau_i:\Sigma\to \Sigma$,
$i\in \{0,1,\ldots,n\} $.
$n$-Quasigroups $f$ and $g$ are called {\em isotopic},
if for some isotopy $\bar \tau=(\tau_0,\tau_1,\ldots ,\tau_n)$ we have
$f(x_1,\dots,x_{n}) \equiv \tau_0^{-1} g(\tau_1 x_{1}, \ldots, \tau_{n} x_{n})$, i.\,e.,
$f\langle x_0,x_1,\dots,x_{n}\rangle \equiv g\langle \tau_0 x_{0}, \tau_1 x_{1}, \ldots, \tau_{n} x_{n}\rangle$.
}
\end{definition}

\begin{definition}\label{def:pr}
{
 An $n$-qua\-si\-group $f$ is termed {\em permutably reducible}
 (in {\rm \cite{PotKro:asymp}}, the term ``decomposable'' was used)
if there exist $m\in\{2,\ldots,n-1\}$,
an $(n-m+1)$-qua\-si\-group $h$,
an $m$-qua\-si\-group  $g$, and a permutation $\sigma: \{1,\ldots,n\} \to \{1,\ldots,n\}$
such that
$$f(x_1,\ldots,x_{n}) \equiv h(g(x_{\sigma(1)},\ldots, x_{\sigma(m)}),
 x_{\sigma(m+1)},\ldots, x_{\sigma(n)})$$
(i.e., $f$ is a {\em composition} of $h$ and $g$).
For short, we will omit the word ``permutably''
(with the exception of the main statements).
If an $n$-qua\-si\-group is not reducible, then it is {\em irreducible}.
(In particular, all $2$-qua\-si\-groups are irreducible.)
}
\end{definition}

\begin{definition}\label{def:SSL}
{
We say that an $n$-qua\-si\-group $f:\{0,1,2,3\}^n\to \{0,1,2,3\}$ is {\em standardly semilinear} if
$$ f\langle \bar x\rangle \leqq L\langle\bar x\rangle
$$
\quad\mbox{where }
$$
L\langle x_0,\ldots,x_{n} \rangle\triangleq l(x_0) \oplus \dots \oplus l(x_{n}) \oplus 1,\quad
l(0)=l(1)=0, \quad l(2)=l(3)=1
$$
($\leqq$ means ``$\leq$ everywhere'';
$\oplus$ means modulo-$2$ addition), see, e.g., Fig.~\ref{fig:2}.
\begin{figure}[ht]
\noindent\mbox{}\hfill
$\setlength\arraycolsep{1mm}
f(x,y):
\begin{array}{|c|c|c|c|}
\hline
  {\mathbf 0} & 1 & 2 & 3 \\
\hline
  1 & 0 & 3 & 2 \\
\hline
  3 & 2 & 0 & 1 \\
\hline
  2 & 3 & 1 & 0 \\ \hline
\end{array}\qquad\qquad
f\langle z,x,y\rangle:
\raisebox{-7mm}%
{
\unitlength 0.5mm
\begin{picture}(38.0,39.2)
\special{em:linewidth 0.3pt}
\definecolor{ogrey}{rgb}{0.5,0.5,0.5}
\color{white}\put(30.0,04.2){\circle*{4.00}}\color{black}\put(30.0,04.2){\circle{4.00}}
\color{white}\put(23.0,04.2){\circle*{4.00}}\color{black}\put(23.0,04.2){\circle{4.00}}
\color{white}\put(16.0,04.2){\circle*{4.00}}\color{black}\put(16.0,04.2){\circle{4.00}}
\color{black}\put(09.0,04.2){\circle*{4.00}}\color{black}\put(09.0,04.2){\circle{4.00}}
\color{white}\put(27.0,02.8){\circle*{4.00}}\color{black}\put(27.0,02.8){\circle{4.00}}
\color{white}\put(20.0,02.8){\circle*{4.00}}\color{black}\put(20.0,02.8){\circle{4.00}}
\color{ogrey}\put(13.0,02.8){\circle*{4.00}}\color{black}\put(13.0,02.8){\circle{4.00}}
\color{white}\put(06.0,02.8){\circle*{4.00}}\color{black}\put(06.0,02.8){\circle{4.00}}
\color{white}\put(24.0,01.4){\circle*{4.00}}\color{black}\put(24.0,01.4){\circle{4.00}}
\color{ogrey}\put(17.0,01.4){\circle*{4.00}}\color{black}\put(17.0,01.4){\circle{4.00}}
\color{white}\put(10.0,01.4){\circle*{4.00}}\color{black}\put(10.0,01.4){\circle{4.00}}
\color{white}\put(03.0,01.4){\circle*{4.00}}\color{black}\put(03.0,01.4){\circle{4.00}}
\color{ogrey}\put(21.0,00.0){\circle*{4.00}}\color{black}\put(21.0,00.0){\circle{4.00}}
\color{white}\put(14.0,00.0){\circle*{4.00}}\color{black}\put(14.0,00.0){\circle{4.00}}
\color{white}\put(07.0,00.0){\circle*{4.00}}\color{black}\put(07.0,00.0){\circle{4.00}}
\color{white}\put(00.0,00.0){\circle*{4.00}}\color{black}\put(00.0,00.0){\circle{4.00}}
\color{white}\put(30.0,13.2){\circle*{4.00}}\color{black}\put(30.0,13.2){\circle{4.00}}
\color{white}\put(23.0,13.2){\circle*{4.00}}\color{black}\put(23.0,13.2){\circle{4.00}}
\color{ogrey}\put(16.0,13.2){\circle*{4.00}}\color{black}\put(16.0,13.2){\circle{4.00}}
\color{white}\put(09.0,13.2){\circle*{4.00}}\color{black}\put(09.0,13.2){\circle{4.00}}
\color{white}\put(27.0,11.8){\circle*{4.00}}\color{black}\put(27.0,11.8){\circle{4.00}}
\color{white}\put(20.0,11.8){\circle*{4.00}}\color{black}\put(20.0,11.8){\circle{4.00}}
\color{white}\put(13.0,11.8){\circle*{4.00}}\color{black}\put(13.0,11.8){\circle{4.00}}
\color{ogrey}\put(06.0,11.8){\circle*{4.00}}\color{black}\put(06.0,11.8){\circle{4.00}}
\color{ogrey}\put(24.0,10.4){\circle*{4.00}}\color{black}\put(24.0,10.4){\circle{4.00}}
\color{white}\put(17.0,10.4){\circle*{4.00}}\color{black}\put(17.0,10.4){\circle{4.00}}
\color{white}\put(10.0,10.4){\circle*{4.00}}\color{black}\put(10.0,10.4){\circle{4.00}}
\color{white}\put(03.0,10.4){\circle*{4.00}}\color{black}\put(03.0,10.4){\circle{4.00}}
\color{white}\put(21.0,09.0){\circle*{4.00}}\color{black}\put(21.0,09.0){\circle{4.00}}
\color{ogrey}\put(14.0,09.0){\circle*{4.00}}\color{black}\put(14.0,09.0){\circle{4.00}}
\color{white}\put(07.0,09.0){\circle*{4.00}}\color{black}\put(07.0,09.0){\circle{4.00}}
\color{white}\put(00.0,09.0){\circle*{4.00}}\color{black}\put(00.0,09.0){\circle{4.00}}
\color{white}\put(30.0,22.2){\circle*{4.00}}\color{black}\put(30.0,22.2){\circle{4.00}}
\color{ogrey}\put(23.0,22.2){\circle*{4.00}}\color{black}\put(23.0,22.2){\circle{4.00}}
\color{white}\put(16.0,22.2){\circle*{4.00}}\color{black}\put(16.0,22.2){\circle{4.00}}
\color{white}\put(09.0,22.2){\circle*{4.00}}\color{black}\put(09.0,22.2){\circle{4.00}}
\color{ogrey}\put(27.0,20.8){\circle*{4.00}}\color{black}\put(27.0,20.8){\circle{4.00}}
\color{white}\put(20.0,20.8){\circle*{4.00}}\color{black}\put(20.0,20.8){\circle{4.00}}
\color{white}\put(13.0,20.8){\circle*{4.00}}\color{black}\put(13.0,20.8){\circle{4.00}}
\color{white}\put(06.0,20.8){\circle*{4.00}}\color{black}\put(06.0,20.8){\circle{4.00}}
\color{white}\put(24.0,19.4){\circle*{4.00}}\color{black}\put(24.0,19.4){\circle{4.00}}
\color{white}\put(17.0,19.4){\circle*{4.00}}\color{black}\put(17.0,19.4){\circle{4.00}}
\color{ogrey}\put(10.0,19.4){\circle*{4.00}}\color{black}\put(10.0,19.4){\circle{4.00}}
\color{white}\put(03.0,19.4){\circle*{4.00}}\color{black}\put(03.0,19.4){\circle{4.00}}
\color{white}\put(21.0,18.0){\circle*{4.00}}\color{black}\put(21.0,18.0){\circle{4.00}}
\color{white}\put(14.0,18.0){\circle*{4.00}}\color{black}\put(14.0,18.0){\circle{4.00}}
\color{white}\put(07.0,18.0){\circle*{4.00}}\color{black}\put(07.0,18.0){\circle{4.00}}
\color{ogrey}\put(00.0,18.0){\circle*{4.00}}\color{black}\put(00.0,18.0){\circle{4.00}}
\color{ogrey}\put(30.0,31.2){\circle*{4.00}}\color{black}\put(30.0,31.2){\circle{4.00}}
\color{white}\put(23.0,31.2){\circle*{4.00}}\color{black}\put(23.0,31.2){\circle{4.00}}
\color{white}\put(16.0,31.2){\circle*{4.00}}\color{black}\put(16.0,31.2){\circle{4.00}}
\color{white}\put(09.0,31.2){\circle*{4.00}}\color{black}\put(09.0,31.2){\circle{4.00}}
\color{white}\put(27.0,29.8){\circle*{4.00}}\color{black}\put(27.0,29.8){\circle{4.00}}
\color{ogrey}\put(20.0,29.8){\circle*{4.00}}\color{black}\put(20.0,29.8){\circle{4.00}}
\color{white}\put(13.0,29.8){\circle*{4.00}}\color{black}\put(13.0,29.8){\circle{4.00}}
\color{white}\put(06.0,29.8){\circle*{4.00}}\color{black}\put(06.0,29.8){\circle{4.00}}
\color{white}\put(24.0,28.4){\circle*{4.00}}\color{black}\put(24.0,28.4){\circle{4.00}}
\color{white}\put(17.0,28.4){\circle*{4.00}}\color{black}\put(17.0,28.4){\circle{4.00}}
\color{white}\put(10.0,28.4){\circle*{4.00}}\color{black}\put(10.0,28.4){\circle{4.00}}
\color{ogrey}\put(03.0,28.4){\circle*{4.00}}\color{black}\put(03.0,28.4){\circle{4.00}}
\color{white}\put(21.0,27.0){\circle*{4.00}}\color{black}\put(21.0,27.0){\circle{4.00}}
\color{white}\put(14.0,27.0){\circle*{4.00}}\color{black}\put(14.0,27.0){\circle{4.00}}
\color{ogrey}\put(07.0,27.0){\circle*{4.00}}\color{black}\put(07.0,27.0){\circle{4.00}}
\color{white}\put(00.0,27.0){\circle*{4.00}}\color{black}\put(00.0,27.0){\circle{4.00}}
\end{picture}
} \leqq
\raisebox{-7mm}%
{
\unitlength 0.5mm
\begin{picture}(38.0,39.2)
\special{em:linewidth 0.3pt}
\definecolor{ogrey}{rgb}{0.5,0.5,0.5}
\color{white}\put(30.0,04.2){\circle*{4.00}}\color{black}\put(30.0,04.2){\circle{4.00}}
\color{white}\put(23.0,04.2){\circle*{4.00}}\color{black}\put(23.0,04.2){\circle{4.00}}
\color{ogrey}\put(16.0,04.2){\circle*{4.00}}\color{black}\put(16.0,04.2){\circle{4.00}}
\color{black}\put(09.0,04.2){\circle*{4.00}}\color{black}\put(09.0,04.2){\circle{4.00}}
\color{white}\put(27.0,02.8){\circle*{4.00}}\color{black}\put(27.0,02.8){\circle{4.00}}
\color{white}\put(20.0,02.8){\circle*{4.00}}\color{black}\put(20.0,02.8){\circle{4.00}}
\color{ogrey}\put(13.0,02.8){\circle*{4.00}}\color{black}\put(13.0,02.8){\circle{4.00}}
\color{ogrey}\put(06.0,02.8){\circle*{4.00}}\color{black}\put(06.0,02.8){\circle{4.00}}
\color{ogrey}\put(24.0,01.4){\circle*{4.00}}\color{black}\put(24.0,01.4){\circle{4.00}}
\color{ogrey}\put(17.0,01.4){\circle*{4.00}}\color{black}\put(17.0,01.4){\circle{4.00}}
\color{white}\put(10.0,01.4){\circle*{4.00}}\color{black}\put(10.0,01.4){\circle{4.00}}
\color{white}\put(03.0,01.4){\circle*{4.00}}\color{black}\put(03.0,01.4){\circle{4.00}}
\color{ogrey}\put(21.0,00.0){\circle*{4.00}}\color{black}\put(21.0,00.0){\circle{4.00}}
\color{ogrey}\put(14.0,00.0){\circle*{4.00}}\color{black}\put(14.0,00.0){\circle{4.00}}
\color{white}\put(07.0,00.0){\circle*{4.00}}\color{black}\put(07.0,00.0){\circle{4.00}}
\color{white}\put(00.0,00.0){\circle*{4.00}}\color{black}\put(00.0,00.0){\circle{4.00}}
\color{white}\put(30.0,13.2){\circle*{4.00}}\color{black}\put(30.0,13.2){\circle{4.00}}
\color{white}\put(23.0,13.2){\circle*{4.00}}\color{black}\put(23.0,13.2){\circle{4.00}}
\color{ogrey}\put(16.0,13.2){\circle*{4.00}}\color{black}\put(16.0,13.2){\circle{4.00}}
\color{ogrey}\put(09.0,13.2){\circle*{4.00}}\color{black}\put(09.0,13.2){\circle{4.00}}
\color{white}\put(27.0,11.8){\circle*{4.00}}\color{black}\put(27.0,11.8){\circle{4.00}}
\color{white}\put(20.0,11.8){\circle*{4.00}}\color{black}\put(20.0,11.8){\circle{4.00}}
\color{ogrey}\put(13.0,11.8){\circle*{4.00}}\color{black}\put(13.0,11.8){\circle{4.00}}
\color{ogrey}\put(06.0,11.8){\circle*{4.00}}\color{black}\put(06.0,11.8){\circle{4.00}}
\color{ogrey}\put(24.0,10.4){\circle*{4.00}}\color{black}\put(24.0,10.4){\circle{4.00}}
\color{ogrey}\put(17.0,10.4){\circle*{4.00}}\color{black}\put(17.0,10.4){\circle{4.00}}
\color{white}\put(10.0,10.4){\circle*{4.00}}\color{black}\put(10.0,10.4){\circle{4.00}}
\color{white}\put(03.0,10.4){\circle*{4.00}}\color{black}\put(03.0,10.4){\circle{4.00}}
\color{ogrey}\put(21.0,09.0){\circle*{4.00}}\color{black}\put(21.0,09.0){\circle{4.00}}
\color{ogrey}\put(14.0,09.0){\circle*{4.00}}\color{black}\put(14.0,09.0){\circle{4.00}}
\color{white}\put(07.0,09.0){\circle*{4.00}}\color{black}\put(07.0,09.0){\circle{4.00}}
\color{white}\put(00.0,09.0){\circle*{4.00}}\color{black}\put(00.0,09.0){\circle{4.00}}
\color{ogrey}\put(30.0,22.2){\circle*{4.00}}\color{black}\put(30.0,22.2){\circle{4.00}}
\color{ogrey}\put(23.0,22.2){\circle*{4.00}}\color{black}\put(23.0,22.2){\circle{4.00}}
\color{white}\put(16.0,22.2){\circle*{4.00}}\color{black}\put(16.0,22.2){\circle{4.00}}
\color{white}\put(09.0,22.2){\circle*{4.00}}\color{black}\put(09.0,22.2){\circle{4.00}}
\color{ogrey}\put(27.0,20.8){\circle*{4.00}}\color{black}\put(27.0,20.8){\circle{4.00}}
\color{ogrey}\put(20.0,20.8){\circle*{4.00}}\color{black}\put(20.0,20.8){\circle{4.00}}
\color{white}\put(13.0,20.8){\circle*{4.00}}\color{black}\put(13.0,20.8){\circle{4.00}}
\color{white}\put(06.0,20.8){\circle*{4.00}}\color{black}\put(06.0,20.8){\circle{4.00}}
\color{white}\put(24.0,19.4){\circle*{4.00}}\color{black}\put(24.0,19.4){\circle{4.00}}
\color{white}\put(17.0,19.4){\circle*{4.00}}\color{black}\put(17.0,19.4){\circle{4.00}}
\color{ogrey}\put(10.0,19.4){\circle*{4.00}}\color{black}\put(10.0,19.4){\circle{4.00}}
\color{ogrey}\put(03.0,19.4){\circle*{4.00}}\color{black}\put(03.0,19.4){\circle{4.00}}
\color{white}\put(21.0,18.0){\circle*{4.00}}\color{black}\put(21.0,18.0){\circle{4.00}}
\color{white}\put(14.0,18.0){\circle*{4.00}}\color{black}\put(14.0,18.0){\circle{4.00}}
\color{ogrey}\put(07.0,18.0){\circle*{4.00}}\color{black}\put(07.0,18.0){\circle{4.00}}
\color{ogrey}\put(00.0,18.0){\circle*{4.00}}\color{black}\put(00.0,18.0){\circle{4.00}}
\color{ogrey}\put(30.0,31.2){\circle*{4.00}}\color{black}\put(30.0,31.2){\circle{4.00}}
\color{ogrey}\put(23.0,31.2){\circle*{4.00}}\color{black}\put(23.0,31.2){\circle{4.00}}
\color{white}\put(16.0,31.2){\circle*{4.00}}\color{black}\put(16.0,31.2){\circle{4.00}}
\color{white}\put(09.0,31.2){\circle*{4.00}}\color{black}\put(09.0,31.2){\circle{4.00}}
\color{ogrey}\put(27.0,29.8){\circle*{4.00}}\color{black}\put(27.0,29.8){\circle{4.00}}
\color{ogrey}\put(20.0,29.8){\circle*{4.00}}\color{black}\put(20.0,29.8){\circle{4.00}}
\color{white}\put(13.0,29.8){\circle*{4.00}}\color{black}\put(13.0,29.8){\circle{4.00}}
\color{white}\put(06.0,29.8){\circle*{4.00}}\color{black}\put(06.0,29.8){\circle{4.00}}
\color{white}\put(24.0,28.4){\circle*{4.00}}\color{black}\put(24.0,28.4){\circle{4.00}}
\color{white}\put(17.0,28.4){\circle*{4.00}}\color{black}\put(17.0,28.4){\circle{4.00}}
\color{ogrey}\put(10.0,28.4){\circle*{4.00}}\color{black}\put(10.0,28.4){\circle{4.00}}
\color{ogrey}\put(03.0,28.4){\circle*{4.00}}\color{black}\put(03.0,28.4){\circle{4.00}}
\color{white}\put(21.0,27.0){\circle*{4.00}}\color{black}\put(21.0,27.0){\circle{4.00}}
\color{white}\put(14.0,27.0){\circle*{4.00}}\color{black}\put(14.0,27.0){\circle{4.00}}
\color{ogrey}\put(07.0,27.0){\circle*{4.00}}\color{black}\put(07.0,27.0){\circle{4.00}}
\color{ogrey}\put(00.0,27.0){\circle*{4.00}}\color{black}\put(00.0,27.0){\circle{4.00}}
\end{picture}
}
$
\hfill\mbox{}
\caption{A standardly semilinear $2$-qua\-si\-group and the function $L\langle\cdot\rangle$}
\label{fig:2}
\end{figure}

An $n$-qua\-si\-group of order $4$ is called {\em semilinear} if
it is isotopic to some standardly semilinear $n$-qua\-si\-group.
}
\end{definition}

The set of standardly semilinear $n$-qua\-si\-groups has a simple characterization:

\begin{proposition}
The following relation is a bijection between the standardly semilinear $n$-qua\-si\-groups
$f$ and the Boolean functions $\lambda:\{0,1\}^n\to \{0,1\}$:
\begin{equation}\label{eq:fB}
\mbox{}~~~~~~
 f\langle x_0,x_1,\ldots,x_n\rangle\equiv
  L\langle x_0,x_1,\ldots,x_n\rangle\cdot
  \big(x_0\oplus x_1\oplus \ldots\oplus x_1\oplus \lambda(l(x_1),\ldots,l(x_n))\big).
\end{equation}
\end{proposition}
\begin{proof}
  Consider an standardly semilinear $n$-qua\-si\-group.
  Consider a set $H$ consisting of $2^{n+2}$ points of $\{0,1,2,3\}^{n+1}$
  with fixed values $l(x_1),\ldots,l(x_n)$.

%
  The number of $1$s of $f\langle\cdot\rangle$ in $H$ is $2^n$
  (indeed, by the definition of an $n$-qua\-si\-group,
  every $n$-tuple $x_1,\ldots,x_n$ corresponds to exactly one $1$).
  Moreover, since $f$ is standardly semilinear,  all these $1$s belong to
  $H_1 \triangleq (\bar x\in H\mid L\langle\bar x \rangle=1)$.
  Since there are no two $1$s that differ in only one coordinate,
  all these $1$s simultaneously have  either even
  or odd coordinate sum. In the even case
  define $\lambda(l(x_1),\ldots,l(x_n))=1$; in the odd case, $=0$.
  Then (\ref{eq:fB}) is automatically true.%
\qquad\end{proof}

So, the number of the standardly
semilinear $n$-qua\-si\-groups is $2^{2^n}$.
Multiplying by the number $3^{n+1}2$ of different functions isotopic to $L$,
we obtain an approximate number of the
semilinear $n$-qua\-si\-groups. The exact number is
$3^{n+1}2^{2^n+1}-8\cdot 6^n$ {\rm \cite[Theorem~1]{PotKro:asymp}},
where $-8\cdot 6^n$ is explained by the fact that affine Boolean functions
(and only affine, i.e., of type
$\lambda(z_1,\ldots,z_n) = b_0 \oplus b_1 z_1 \oplus \ldots \oplus b_n z_n$, $b_i\in \{0,1\}$ )
corresponds to $n$-qua\-si\-groups majorized by more than one isotope of $L$.

In the rest of the paper, unless otherwise stated,
we consider only order-$4$ $n$-qua\-si\-groups over $\Sigma=\{0,1,2,3\}$.

\section{Main result}\label{s:2}

The main result is the following theorem.
\begin{theorem} \label{th:main}
Every $n$-qua\-si\-group of order $4$ is permutably reducible or semilinear.
\end{theorem}

The basic characteristic of an $n$-qua\-si\-group $f$,
which divides our proof into four subcases,
is the maximum arity of its irreducible retract.
Denote this value by $\kappa(f)$; then, $2\leq\kappa(f)<n$.
The line of reasoning in the proof of Theorem~\ref{th:main} is inductive,
so we can assume that the irreducible retracts are semilinear.

\begin{lemma}[Case $\kappa=n-1$, {\rm \cite[Lemma~4]{PotKro:asymp}}] \label{pro:3b}
If an $n$-qua\-si\-group $f$ of order $4$ has a semilinear $(n-1)$-retract,
then it is permutably reducible or semilinear.
\end{lemma}

\begin{lemma}[Case $2<\kappa\leq n-3$, {\rm \cite{Kro:n-3}}] \label{pro:2b}
Let $f$ be an $n$-qua\-si\-group of arbitrary order and $\kappa(f)\in\{3,\ldots,n-3\}$.
Then $f$ is permutably reducible.
\end{lemma}
In \cite{Kro:n-2}, an example of irreducible $n$-qua\-si\-group of order $4$
whose $(n-1)$-retracts are all reducible is constructed for every even $n\geq 4$.
So, the assumption of Lemma~\ref{pro:2b} cannot be extended
to the case $\kappa = n-2$.
Nevertheless, in Section~\ref{sect:4b} we will prove the following:

\begin{lemma}[Case $\kappa=n-2$] \label{pro:4b}
Let $n\geq 5$.
If an $n$-qua\-si\-group $f$ of order $4$ has a semilinear permutably irreducible $(n-2)$-retract
and all the $(n-1)$-retracts are permutably reducible,
then $f$ is permutably reducible or semilinear.
\end{lemma}

The last case, announced in \cite{Pot:2006:thes}, will be proved in Section~\ref{sect:1b}:
\begin{lemma}[Case $\kappa=2$] \label{pro:1b}
Let $n\geq 5$; let $f$ be an $n$-qua\-si\-group of order $4$,
 and let all its $k$-retracts with $2<k<n$ be permutably reducible.
 Then $f$ is permutably reducible.
\end{lemma}

{\em Proof of Theorem~\rm\ref{th:main}}.
The validity of the Theorem for $n\leq 4$ (and even for $n\leq 5$) is proved by exhaustion.
Assume, by induction, that all $m$-qua\-si\-groups of order $4$ with $m<n$ are reducible or semilinear.
Consider an $n$-qua\-si\-group $f$ of order $4$.
It has an irreducible $\kappa(f)$-retract, which is semilinear by inductive assumption.
Depending on the value of $\kappa(f)= 2, 3\ldots n{-}3,n{-}2,n{-}1$,
the statement of Theorem follows from one of Lemmas
\ref {pro:1b}, \ref {pro:2b}, \ref {pro:4b}, \ref {pro:3b}.\qquad
\endproof


\section{Proof of Lemma~\ref{pro:1b}}\label{sect:1b}\label{s:3}

We will prove a stronger variant (Lemma~\ref{pro:1b-alt}) of the statement.
It uses the following concept:
\begin{definition}
{
An $n$-qua\-si\-group $f$ is called  {\em completely reducible}
if it is permutably reducible and all its principal retracts
of arity more than $2$ are permutably reducible (equivalently,
$f$ can be represented as a composition of $n-1$ binary quasigroups, e.g.
$f(x_1,x_2,x_3,x_4,x_5)=f_1(f_2(x_1,x_3),f_3(f_4(x_2,x_5),x_4))$\,).
}
\end{definition}

\begin{lemma}\label{pro:1b-alt}
Let all the principal $3$- and $4$-retracts of an
$n$-qua\-si\-group $f$ of order $4$ ($n\geq 5$) be permutably reducible. Then
$f$ is completely reducible.
\end{lemma}

\begin{definition}
{
An $n$-qua\-si\-group $f$ is called {\em normalized}
if for all $i\in \{1,\ldots,n\}$ and $a\in \Sigma$ it is true that $f(\bar 0^{\laaa i\raaa}\lbbb a\rbbb) =a$.
}
\end{definition}

Denote by $\Gamma$ the set of all (four) normalized binary
qua\-si\-groups of order $4$. It is straightforward that the
operations from $\Gamma$ are associative and commutative (they are
isomorphic to the additive groups $Z_2^2$ and $Z_4$), and we will use the
 form $a \star b$ instead of $\star(a,b)$ to write the result of $\star \in \Gamma$.

 Let $K_n = \langle V(K_n),E(K_n)\rangle$ be the complete graph with $n$
 vertices associated with the arguments $x_1,\ldots,x_n$ of an $n$-ary operation.
 For the edges, we will use the short notation like $x_i x_j$.
 For any normalized $n$-qua\-si\-group $f$
we define the edge coloring $\mu_f:E(K_n)\to\Gamma$
by the following way: the color
$\mu_f(x_i x_j)$ of an edge $x_i x_j \in E(K_n)$ is defined as
the binary operation $\star$ such that
$f(\bar 0^{\laaa i,j\raaa} \lbbb x_i,x_j\rbbb)\equiv x_i\star x_j$.
%

 \begin{proposition} \label{pro:21}
 Let $f$ be an $n$-qua\-si\-group, $n\geq 5$, and
 let all $3$- and $4$-retracts of $f$ be reducible.
 Then the coloring $\mu_f$ of $K_n$ satisfies the following:
\begin{itemize}
 \item[\rm (A)] Every triangle is colored by at most $2$ colors.
 \item[\rm (B)] If a tetrahedron is colored by $2$ colors with $3$ edges of each color,
  then it includes a one-color triangle; i.\,e., the following fragment is forbidden:\\[2mm]
\hspace{0em}\raisebox{-6mm}{
\def\dick{2.0pt}
\def\duenn{0.4pt}
\unitlength 0.3mm
\def\emline#1#2#3#4#5#6{%
\put(#1,#2){\special{em:moveto}}%
\put(#4,#5){\special{em:lineto}}%
}
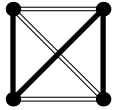
\begin{picture}(40.00,40.00)
\color[rgb]{0.000,0.000,0.000}
\special{em:linewidth \duenn}
\emline{0.00}{ 0.80}{1}{40.00}{  0.80}{2}
\emline{0.00}{-0.80}{1}{40.00}{ -0.80}{2}
\emline{0.00}{40.80}{1}{40.00}{ 40.80}{2}
\emline{0.00}{39.20}{1}{40.00}{ 39.20}{2}
\emline{0.60}{40.60}{1}{40.60}{  0.60}{2}
\emline{-0.6}{39.40}{1}{39.40}{ -0.60}{2}
\color[rgb]{0.000,0.000,0.000}
\special{em:linewidth \dick}
\emline{40.00}{ 40.00}{1}{40.00}{  0.00}{2}
\emline{ 0.00}{  0.00}{1}{ 0.00}{ 40.00}{2}
\emline{ 0.00}{  0.00}{1}{40.00}{ 40.00}{2}
\special{em:linewidth 5.5pt}
\emline{ 0.00}{ 40.00}{1}{ 0.00}{ 40.00}{2}
\emline{40.00}{ 40.00}{1}{40.00}{ 40.00}{2}
\emline{ 0.00}{  0.00}{1}{ 0.00}{  0.00}{2}
\emline{40.00}{  0.00}{1}{40.00}{  0.00}{2}
\end{picture}
}
\end{itemize}
  \end{proposition}
  \vspace{2mm}

  \begin{proof}
  Every $3$-retract of $f$ is a composition of two binary operations, which yields (A).

 Consider the $4$-retract $f_4$ of $f$ that
corresponds to some four vertices (the other variables are fixed by $0$).
Since it is reducible and normalized, it can be represented as a composition of some
normalized $2$-qua\-si\-group and $3$-qua\-si\-group.
The $3$-qua\-si\-group is a $3$-retract of $f$ and,
in its turn, can be represented as a composition of
two normalized $2$-qua\-si\-groups. So, $f_4(x,y,z,u)$, up to permutation of arguments,
has the form $(x\star y)\circ(u\diamond v)$ or $x\circ (y\star (u\diamond v))$
for some $\star,\circ,\diamond\in\Gamma$.
As follows from the hypothesis of (B), two of these three operations coincide.
Thus, there are only four types of decomposition of $f_4$:
  $(x\star y)\circ(u\star v)$, $x\star(y\circ u\circ v)$, $x\star y\star (u\circ v)$,
   $x\circ (y\star (u\circ v))$.
   In any case, (B) holds.\qquad
   \end{proof}

\begin{proposition} \label{cor:21}
Assume that an edge coloring $\mu$ of $K_n$ satisfies {\rm (A)} and {\rm (B)}.
Then for each pairwise different $a, b, c, d \in V(K_n)$  the condition
$\mu(ab)=\mu(ac)\neq\mu(bc)=\mu(bd)\neq\mu(cd)$ implies
$\mu(ad)=\mu(ab)$. I.\,e.,
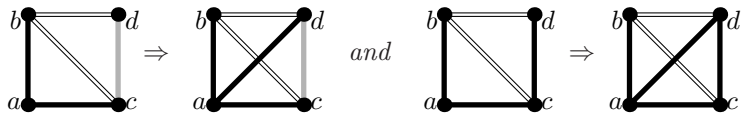
$$
\def\dick{2.0pt}
\def\duenn{0.4pt}
\unitlength 0.3mm
\def\emline#1#2#3#4#5#6{%
\put(#1,#2){\special{em:moveto}}%
\put(#4,#5){\special{em:lineto}}}
\def\pictSample{
\color[rgb]{0.000,0.000,0.000}
\emline{ 0.00}{-20.00}{1}{ 0.00}{ 20.00}{2}
\emline{ 0.00}{-20.00}{1}{40.00}{-20.00}{2}
\special{em:linewidth 5.5pt}
\emline{ 0.00}{ 20.00}{1}{ 0.00}{ 20.00}{2}
\emline{40.00}{ 20.00}{1}{40.00}{ 20.00}{2}
\emline{ 0.00}{-20.00}{1}{ 0.00}{-20.00}{2}
\emline{40.00}{-20.00}{1}{40.00}{-20.00}{2}
\special{em:linewidth \duenn}
\emline{0.00}{20.80}{1}{40.00}{ 20.80}{2}
\emline{0.00}{19.20}{1}{40.00}{ 19.20}{2}
\emline{0.60}{20.60}{1}{40.60}{-19.40}{2}
\emline{-0.6}{19.40}{1}{39.40}{-20.60}{2}
\put(-3.00,22.00){\makebox(0,0)[tr]{$b$}}
\put(-3.00,-22.00){\makebox(0,0)[br]{$a$}}
\put(43.00,22.00){\makebox(0,0)[tl]{$d$}}
\put(43.00,-22.00){\makebox(0,0)[bl]{$c$}} }
\begin{picture}(40.00,18.00)
\color[rgb]{0.000,0.000,0.000}
\special{em:linewidth \dick}
\color[rgb]{0.702,0.702,0.702}
\emline{40.00}{20.00}{1}{40.00}{-20.00}{2}
\pictSample
\end{picture}
\quad\Rightarrow\quad
\begin{picture}(40.00,18.00)
\special{em:linewidth \dick}
\color[rgb]{0.702,0.702,0.702}
\emline{40.00}{20.00}{1}{40.00}{-20.00}{2}
\pictSample
\special{em:linewidth \dick}
\emline{0.00}{-20.00}{1}{40.00}{20.00}{2}
\end{picture}
\qquad\mbox{and}\qquad
\begin{picture}(40.00,18.00)
\color[rgb]{0.000,0.000,0.000}
\special{em:linewidth \dick}
\emline{40.00}{20.00}{1}{40.00}{-20.00}{2}
\pictSample
\special{em:linewidth \dick}
\end{picture}
\quad\Rightarrow\quad
\begin{picture}(40.00,18.00)
\color[rgb]{0.000,0.000,0.000}
\special{em:linewidth \dick}
\emline{40.00}{20.00}{1}{40.00}{-20.00}{2}
\pictSample
\special{em:linewidth \dick}
\emline{0.00}{-20.00}{1}{40.00}{20.00}{2}
\end{picture}
$$\vspace{0mm}
\end{proposition}
\begin{proof}
Obviously, any other variant for $\mu(ad)$ contradicts  (A) or (B).\qquad
\end{proof}

The following proposition is easy to check.
\begin{proposition} \label{pro:24}
 Let $f$ and $g$ be reducible $3$-qua\-si\-groups of order $4$, and
 let
 $ f(x,y,0)\equiv g(x,y,0)$,
 $ f(x,0,z)\equiv g(x,0,z)$, and
 $ f(0,y,z)\equiv g(0,y,z)$.
 Then $f\equiv g$.
\end{proposition}
\begin{remark} \label{n:1}
{\rm
Indeed, Proposition~\ref{pro:24} holds for every order with the
extra condition that $f$ is a composition of two different or
associative $2$-qua\-si\-groups. The similar statement for
$n$-qua\-si\-groups with $n>3$ holds for an arbitrary order without
extra conditions {\rm \cite[Theorem~1]{KPS:ir}}: if two reducible $n$-qua\-si\-groups
coincide on every $n$-tuples with one zero, then they are identical.
}
\end{remark}
\begin{corollary} \label{cor:24}
 Let $f$ and $g$ be $n$-qua\-si\-groups of order $4$ ($n\geq 3$)
 whose principal $3$-retracts are all reducible.
 Assume that
 $f(\bar 0^{\laaa i,j\raaa}\lbbb y,z\rbbb)\equiv g(\bar 0^{\laaa i,j\raaa}\laaa y,z\raaa)$
 for every $i,j\in \{1,\ldots,n\}$ and $y,z\in\Sigma$.
 Then $f\equiv g$.
\end{corollary}
\begin{proof} The equality $f(\bar x)\equiv g(\bar x)$ is proved by
induction on the number of non-zero elements in $\bar x$, using the
reducibility of $3$-retracts and Proposition~\ref{pro:24}.

For example, to prove that $f(1,2,3,2,1,\bar 0) = g(1,2,3,2,1,\bar 0)$
we can consider the $3$-retracts
$f^3(x,y,z)\triangleq f(1,2,x,y,z,\bar 0)$
and
$g^3(x,y,z)\triangleq g(1,2,x,y,z,\bar 0)$.
By the induction assumption $f^3$ and $g^3$ meet the hypothesis of Proposition~\ref{pro:24}.
Thus, $f^3 \equiv g^3$, and, in particular, $f^3(3,2,1)=g^3(3,2,1)$.%
\qquad
\end{proof}

The following proposition is the key statement in the proof of Lemma~\ref{pro:1b-alt}.
\begin{proposition} \label{pro:22}
 Assume that an edge coloring $\mu:E(K_n)\to\Gamma$ of the graph $K_n$ meets {\rm (A)} and {\rm (B)}.
 Then there exists a completely reducible $n$-qua\-si\-group $f$
 such that $\mu_f= \mu$.
\end{proposition}

Before proving Proposition~\ref{pro:22} by induction,
we consider one auxiliary statement,
which will be used in the induction step.
We say that an edge $xy\in E(K_n)$ is \emph{inner} with respect to some
edge coloring $\mu$ of $K_n$ if
  for any $z\in V(K_n)\setminus \{x,y\}$ it is true that $\mu(xz)=\mu(yz)$.

\begin{proposition} \label{pro:*}
 Assume that an edge coloring $\mu$ of $K_n$ meets {\rm (A)} and {\rm (B)}.
 Then $K_n$ contains an inner edge.
\end{proposition}
\begin{proof}
Consider an
arbitrary sequence of edges
$e_1,e_2,...,e_k$ that satisfies the following:
\begin{itemize}
\item[(C)] for every $j\in \{1,...,k-1\}$ the edges $e_{j}$ and $e_{j+1}$ are
adjacent and
$\mu(e_{j})=\mu(e_{j}\triangle e_{j+1})\neq \mu(e_{j+1})$
(where $\triangle$ means the symmetrical difference between two sets).
\end{itemize}
Denote by $a_j$ the element from $e_j \setminus e_{j+1}$.
\begin{figure}[ht]
$$
\def\dick{4pt}
\def\normal{1.5pt}
\def\duenn{0.3pt}
\unitlength 0.6mm
\def\emline#1#2#3#4#5#6{\put(#1,#2){\special{em:moveto}}\put(#4,#5){\special{em:lineto}}}
\raisebox{-12mm}{
\begin{picture}(66,48.00)
\special{em:linewidth \dick}
\emline{50.00}{ 0.00}{1}{62.00}{20.00}{2}
\color[rgb]{0.7,0.7,0.7}
\emline{50.00}{40.00}{1}{50.00}{  0.00}{2}
\special{em:linewidth \normal}
\emline{50.00}{40.00}{1}{62.00}{ 20.00}{2}
\color[rgb]{0.000,0.000,0.000}
\special{em:linewidth \duenn}
\emline{30.56}{ -0.44}{1}{50.56}{39.56}{2}
\emline{29.44}{  0.44}{1}{49.44}{40.44}{2}
\emline{30.00}{  0.44}{1}{50.00}{ 0.44}{2}
\emline{30.00}{ -0.44}{1}{50.00}{-0.44}{2}
\special{em:linewidth \dick}
\emline{30.00}{40.00}{1}{30.00}{  0.00}{2}
\special{em:linewidth \normal}
\emline{30.00}{40.00}{1}{50.00}{ 40.00}{2}
\color[rgb]{0.7,0.7,0.7}
\emline{18.00}{22.00}{1}{30.00}{ 40.00}{2}
\special{em:linewidth \dick}
\emline{18.00}{20.00}{1}{30.00}{  0.00}{2}
\color[rgb]{0.000,0.000,0.000}
\emline{ 0.00}{20.00}{1}{18.00}{ 22.00}{2}
\special{em:linewidth \normal}
\emline{ 0.00}{20.00}{1}{30.00}{  0.00}{2}
\special{em:linewidth \duenn}
\put( 0.00,20.00){\circle*{ 4.00}}
\put(18.11,22.00){\circle*{ 4.00}}
\put(30.11,40.00){\circle*{ 4.00}}
\put(30.11, 0.00){\circle*{ 4.00}}
\put(50.11,40.00){\circle*{ 4.00}}
\put(50.11, 0.00){\circle*{ 4.00}}
\put(61.78,20.00){\circle*{ 4.00}}
\put(23.00,12.00){\makebox(0,0)[bl]{$e_2$}}
\put(31.10,20.00){\makebox(0,0)[cl]{$e_3$}}
\put(37.60,14.40){\makebox(0,0)[tl]{$e_4$}}
\put(49.40,18.60){\makebox(0,0)[cr]{$e_5$}}
\put(54.60, 6.40){\makebox(0,0)[tl]{$e_6$}}
\put( 7.60,21.40){\makebox(0,0)[bc]{$e_1$}}
\end{picture}}
\quad\Rightarrow\quad \raisebox{-12mm}{
\begin{picture}(66,48.00)
\special{em:linewidth \dick}
\emline{50.00}{ 0.00}{1}{62.00}{20.00}{2}
\color[rgb]{0.7,0.7,0.7}
\emline{50.00}{40.00}{1}{50.00}{  0.00}{2}
\special{em:linewidth \normal}
\emline{50.00}{40.00}{1}{62.00}{ 20.00}{2}
\color[rgb]{1.000,1.000,1.000}
\special{em:linewidth \dick}
\emline{30.00}{ 0.00}{1}{50.00}{40.00}{2}
\special{em:linewidth\normal}
\emline{30.00}{  0.00}{1}{50.00}{ 0.00}{2}
\emline{30.00}{0.00}{1}{61.66}{20.00}{2}
\color[rgb]{0.000,0.000,0.000}
\special{em:linewidth \duenn}
\emline{30.56}{-0.44}{1}{50.56}{39.56}{2}
\emline{29.44}{0.44}{1}{49.44}{40.44}{2}
\emline{30.00}{0.44}{1}{50.00}{ 0.44}{2}
\emline{30.00}{-0.44}{1}{50.00}{-0.44}{2}
\emline{30.33}{-0.33}{1}{62.00}{19.67}{2}
\emline{29.67}{0.33}{1}{61.33}{20.33}{2}
\special{em:linewidth \dick}
\emline{30.00}{40.00}{1}{30.00}{ 0.00}{2}
\special{em:linewidth \normal}
\emline{30.00}{40.00}{1}{50.00}{40.00}{2}
\emline{30.00}{40.00}{1}{42.80}{32.00}{2}
\emline{30.00}{40.00}{1}{38.00}{24.00}{2}
\color[rgb]{0.7,0.7,0.7}
\special{em:linewidth \normal}
\emline{18.00}{22.00}{1}{30.00}{40.00}{2}
\emline{18.00}{22.00}{1}{28.00}{27.66}{2}
\emline{18.00}{22.00}{1}{28.00}{21.55}{2}
\emline{18.00}{22.00}{1}{28.00}{15.13}{2}
\special{em:linewidth \dick}
\emline{18.00}{20.00}{1}{30.00}{0.00}{2}
\color[rgb]{0.000,0.000,0.000}
\emline{0.00}{20.00}{1}{18.00}{22.00}{2}
\special{em:linewidth \normal}
\emline{0.00}{20.00}{1}{30.00}{ 0.00}{2}
\emline{0.00}{20.00}{1}{27.33}{38.33}{2}
\emline{0.00}{20.00}{1}{20.67}{30.00}{2}
\emline{0.00}{20.00}{1}{20.33}{10.00}{2}
\emline{0.00}{20.00}{1}{15.33}{20.00}{2}
\special{em:linewidth \duenn}
\put( 0.00,20.00){\circle*{ 4.00}}
\put(18.11,22.00){\circle*{ 4.00}}
\put(30.11,40.00){\circle*{ 4.00}}
\put(30.11, 0.00){\circle*{ 4.00}}
\put(50.11,40.00){\circle*{ 4.00}}
\put(50.11, 0.00){\circle*{ 4.00}}
\put(61.78,20.00){\circle*{ 4.00}}
\put( 0.00,22.40){\makebox(0,0)[bc]{$a_1$}}
\put(17.00,23.00){\makebox(0,0)[br]{$a_2$}}
\put(27.40,39.00){\makebox(0,0)[br]{$a_3$}}
\put(52.60,38.40){\makebox(0,0)[bl]{$a_5$}}
\put(27.40, 0.40){\makebox(0,0)[tr]{$a_4$}}
\end{picture}}
$$
\caption{An example of a sequence $e_1,\ldots,e_k$ from the proof of Proposition~\ref{pro:22}}
\label{fig:1}
\end{figure}
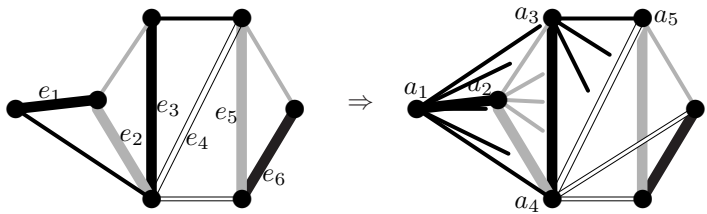

\emph{Claim~(*)}. We claim that for every $i$, $j$, $1\leq i < j\leq k$, and
$d\in e_j$ the vertices $a_i$ and $d$ are different and $\mu(a_i
d)=\mu(e_i)$ (see Fig.~\ref{fig:1}). We will show this by induction on $j-i$.

If $j-i=1$, the claim follows from (C). If $j-i=2$, the claim
follows from Proposition~\ref{cor:21} ($a:=a_i$, $b:=a_{i+1}$).
Assume $j-i>2$. By the inductive assumption,
\begin{eqnarray*}
&&\mu(e_i)=\mu(a_i a_{i+1})=\mu(a_i a_{i+2}) \\
&\neq&
\mu(e_{i+1})=\mu(a_{i+1} a_{i+2})=\mu(a_{i+1} d)\\
&\neq& \mu(e_{i+2})=\mu(a_{i+2} d).
\end{eqnarray*}
Consequently, $a_i\neq d$, and, by Proposition~\ref{cor:21} ($a:=a_i$, $b:=a_{i+1}$, $c:=a_{i+2}$),
we have $\mu(a_i d)=\mu(e_i)$. Claim (*) is proved.

So, all $a_1$, \ldots, $a_{k-1}$ are mutually different, and thus
there exists a maximum sequence $e_1, e_2,\ldots,e_k$ satisfying (C).
Then, its maximality
and (*) imply that the edge $e_k$ is inner.%
\qquad
\end{proof}

{\em Proof of Proposition~\rm\ref{pro:22}}.
We will proceed by induction on $n$.
If $n=3$, then the statement is trivially true.
Assume that Proposition~\ref{pro:22} holds for $(n-1)$-qua\-si\-groups.
Consider an inner edge $e\in E(K_n)$.
Without loss of generality we can assume $e=x_{n-1}x_n$, $\mu(e)=\star$.
Denote by $\mu^{n-1}$ the restriction of the coloring
$\mu$ on $K_{n-1}\subset K_n$. By the inductive assumption, there exists
a completely reducible $(n-1)$-qua\-si\-group $g$ such that
$\mu_g=\mu^{n-1}$. Then $f(x_1,\dots,x_n) \triangleq
g(x_1,x_2,\dots, x_{n-2},x_{n-1}\star x_n)$ is a desired
$n$-qua\-si\-group
(indeed, $\mu_f(x_{n-1}x_n)=\star=\mu(x_{n-1}x_n)$;
if $i<j<n$, then $\mu_f(x_i x_j)=\mu_g(x_i x_j)=\mu(x_i x_j)$;
if $i<n-1$, then $\mu_f(x_i x_n)=\mu_g(x_i x_{n-1})=\mu(x_i x_{n-1})=\mu(x_i x_{n})$,
where the last equality follows from the innerness of $x_{n-1}x_n$).%
\qquad\endproof

{\em Proof of Lemma~\rm\ref{pro:1b-alt}}.
Let $f$ be an
$n$-qua\-si\-group of order $4$
whose $3$- and $4$-retracts are all reducible.
Without loss of generality we assume that $f$ is
normalized (otherwise, we can normalize it, applying an appropriate
isotopy). Then, by Proposition~\ref{pro:21}, the corresponding
edge coloring $\mu_f$ of the graph $K_n$ satisfies (A) and (B).
By Proposition~\ref{pro:22}, there exists a completely reducible
$n$-qua\-si\-group $g$ with $\mu_g\equiv \mu_f$.
By Corollary~\ref{cor:24}, $f$ and $g$ are identical.%
\qquad\endproof

\section{Proof of Lemma~\ref{pro:4b}%
}\label{sect:4b}\label{s:4}
In the proof, we will use the following three propositions.
The first simple one, on a representation of a reducible $n$-qua\-si\-group with an irreducible $(n-1)$-retract,
holds for an arbitrary order.
\begin{proposition} \label{p:canon}
 Assume that a reducible $n$-qua\-si\-group $D$ ($n\geq 3$) of an arbitrary order
 has an irreducible $(n-1)$-retract
 $F\langle x_0,\ldots,x_{n-1}\rangle \equiv  D\langle x_0,\ldots,x_{n-1},0\rangle$.
Then there are $i\in \{0,\ldots ,n\}$ and a $2$-qua\-si\-group $h$
such that $h(x,0)\equiv x$ and
\begin{equation}\label{eq:canon}
D\langle x_0,\ldots,x_{n}\rangle\equiv F\langle x_0,\ldots, x_{i-1},h(x_i,x_{n}), x_{i+1},\ldots, x_{n-1}\rangle.
\end{equation}
\end{proposition}

\begin{proof}
Since $D$ is reducible, $D\langle x_0,\ldots,x_{n}\rangle$ can be represented as
$H\langle f(\bar x'),\bar x'' \rangle$ where $\bar x'$ and $\bar x''$ are disjoint
groups of variables, each containing at least two variables.
If $x_n$ is grouped with more than one other variable, then fixing $x_n$ gives
a reducible retract, which contradicts the irreducibility of $F$.
So, we conclude that for some $i\in \{0,\ldots n-1\}$ there exists one of the following
two representations of $D$:
\begin{eqnarray}\label{eq:first}
D\langle x_0,\ldots,x_{n}\rangle &\equiv& G\langle g(x_i,x_{n}),\tilde x\rangle,\\
D\langle x_0,\ldots,x_{n}\rangle &\equiv& g\langle G(\tilde x),x_i,x_{n}\rangle\nonumber
\end{eqnarray}
where $\tilde x \triangleq (x_0,\ldots ,x_{i-1},x_{i+1},\ldots ,x_{n-1})$ and
$G$ and $g$ are $(n-1)$- and $2$- qua\-si\-groups.
Moreover, the existence of a representation of the first type implies
the existence of a representation of the second type, and vice versa:
$$
g\langle G(\tilde x),x_i,x_{n}\rangle
\equiv
\cases{ 1& if $G(\tilde x) =   g(x_i,x_{n})$ \cr
        0& if $G(\tilde x)\neq g(x_i,x_{n})$     }
\equiv
G\langle g(x_i,x_{n}),\tilde x\rangle.
$$
So, we can assume that (\ref{eq:first}) holds. Put $\gamma(x_i)\triangleq g(x_i,0)$.
Then, $F\langle x_0,\ldots,x_{n-1}\rangle \equiv G\langle \gamma(x_i), \tilde x \rangle$,
and (\ref{eq:canon}) holds with $h(x_i,x_{n})\triangleq \gamma^{-1}g(x_i,x_{n})$.%
\qquad\end{proof}

In what follows,
permutations $\sigma:\Sigma \to \Sigma$ will be denoted by the value
lists $(\sigma(0),\sigma(1),\sigma(2),\sigma(3))$; denote $Id \triangleq (0{,}1{,}2{,}3)$.

\begin{proposition}[on autotopies of a semilinear $n$-qua\-si\-group] \label{p:aut}
Assume $f$ is a standardly semilinear $n$-qua\-si\-group. Denote by $\pi$ the permutation
$(1{,}0{,}3{,}2)$.
Then for every different $i,j\in \{0,\ldots,n\}$ the following holds:\\
{\rm a)} $f\langle \bar x \rangle \equiv f\langle \bar x^{\laaa i,j \raaa}{\lbbb \pi x_i, \pi x_j \rbbb} \rangle$,
where $\bar x\triangleq (x_0,\ldots ,x_n)$;\\
{\rm b)} if $f\langle \bar x \rangle \equiv f\langle \bar x^{\laaa i,j \raaa}{\lbbb \mu x_i, \nu x_j \rbbb} \rangle$
holds for some other pair of non-identity permutations $(\mu,\nu)\neq(\pi,\pi)$ and $n\geq 3$, then $f$ is reducible.
\end{proposition}

\begin{proof}
a) It is straightforward that
$f\langle \bar x^{\laaa i \raaa}{\lbbb \pi x_i \rbbb}\rangle \equiv L\langle \bar x \rangle-f\langle \bar x \rangle$,
where $L\langle \cdot \rangle$ is from Definition~\ref{def:SSL}.
So, $f\langle \bar x^{\laaa i,j \raaa}{\lbbb \pi x_i, \pi x_j \rbbb}\rangle
\equiv L\langle \bar x \rangle-(L\langle \bar x \rangle-f\langle \bar x \rangle)\equiv f\langle \bar x \rangle$.

b) Without loss of generality assume that $i=1$, $j=2$.
Put
\begin{eqnarray*}
\alpha(x,y) &\triangleq& f(x,y,\bar 0), 
\\
\beta(x,\bar z) &\triangleq& f(x,0,\bar z),\qquad \bar z\triangleq(z_1,\ldots ,z_{n-2}),\\
\gamma(x) &\triangleq& f(x,0,\bar 0). 
\end{eqnarray*}
Assume there exists a pair $(\mu, \nu)$ that satisfies the hypothesis of b).
Then $\alpha(x,y) \equiv \alpha(\mu x,\nu y)$.
It is easy to see 
that the permutation $\nu$ does not have fixed points.
So, $\nu$ is either a cyclic permutation or a involution
($(2{,}3{,}0{,}1)$ or $(3{,}2{,}1{,}0)$) different from $\pi=(1{,}0{,}3{,}2)$.
In any case, we can derive the following:

\emph{Claim~(*)}. For each
$v\in\Sigma$ there exist permutations $\rho_v,\tau_v:\Sigma\to \Sigma$
such that $f(x,y,\bar z)\equiv f(\rho_v x,\tau_v y,\bar z)$ and $\tau_v v=0$
(in other words, the group of permutations $\tau$ admitting $f(x,y,\bar z)\equiv f(\rho x,\tau y,\bar z)$
for some $\rho$ acts transitively on $\Sigma$, i.e., has only one orbit):

\emph{Case}~1. If $\nu$ is a cyclic permutation,
then $v$, $\nu v$, $\nu^2 v$, $\nu^3 v$ are pairwise different;
so, one of the pairs $(Id,Id)$, $(\mu,\nu)$, $(\mu^2,\nu^2)$, $(\mu^3,\nu^3)$
can be chosen as $(\rho_v,\tau_v)$, proving (*).

\emph{Case}~2. If $\nu$ is $(2{,}3{,}0{,}1)$ or $(3{,}2{,}1{,}0)$,
then $v$, $\nu v$, $\pi v$, $\nu \pi v$ are pairwise different,
and $(\rho_v,\tau_v)$ can be chosen from
$(Id,Id)$, $(\mu,\nu)$, $(\pi,\pi)$, $(\mu \pi,\nu \pi)$.

Claim (*) is proved. Then,
\begin{eqnarray*}
f(x,y,\bar z)
&\equiv&
f(\rho_y x,\tau_y y,\bar z)
\equiv
f(\rho_y x,0,\bar z)
\equiv
\beta(\rho_y x,\bar z)
\\&\equiv&
\beta(\gamma^{-1}\alpha(\rho_y x,0),\bar z)
\equiv
\beta(\gamma^{-1}\alpha(\rho_y x,\tau_y y),\bar z)
\equiv
\beta(\gamma^{-1}\alpha(x,y),\bar z)
\end{eqnarray*}
and thus $f$ is reducible provided $n\geq 3$.%
\qquad\end{proof}

The next proposition concerns $2$-qua\-si\-groups of order $4$, and the proof is straightforward.
\begin{proposition} \label{p:2-2}
Let $s$ and $t$ be $2$-qua\-si\-groups.
Denote $s_i(x) \triangleq s(x,i)$ and $t_i(x) \triangleq t(x,i)$.
Let $s_0= t_0= Id$, and let
for every $i$ either $t_i s_i^{-1} = Id$ or $t_i s_i^{-1}  = (1{,}0{,}3{,}2)$. 
Then either $s\equiv t$ or
for some permutation $\phi$ the $2$-qua\-si\-group $s'(x,y)\triangleq s( x,\phi y)$
is standardly semilinear.
\end{proposition}

\begin{proof}
Denote $\pi \triangleq (1{,}0{,}3{,}2)$; observe $\pi=\pi^{-1}$.
Assume that $s\not\equiv t$.
Then there are at least two different elements $i,j\in\{1,2,3\}$ for which $t_i s_i^{-1}  = \pi$.
Denote by $k$ the third element, i.e., $\{i,j,k\}=\{1,2,3\}$.
The permutation $t_i$ has no fixed points, otherwise there is a contradiction with $t_0$;
similarly, $s_i=\pi t_i$ has no fixed points.
So, $t_i$ and, similarly, $t_j$ belong to $\{(2{,}3{,}0{,}1)$, $(2{,}3{,}1{,}0)$, $(3{,}2{,}0{,}1)$, $(3{,}2{,}1{,}0)\}$.
The only variant for $t_k$ is $(1{,}0{,}3{,}2)$.
Then, $s'(x,y)\triangleq s( x,\phi y)$ is standardly semilinear
with $\phi \triangleq (0{,}k{,}i{,}j)$.%
\qquad\end{proof}

{\em Proof of Lemma~\rm\ref{pro:4b}}.
Assume $C$ is an $n$-qua\-si\-group of order $4$.
Assume all the $(n-1)$-retracts of $C$ are reducible
and $C$ has a semilinear irreducible $(n-2)$-retract $E$.
Without loss of generality assume that
$$ E\langle x_0,\ldots,x_{n-2} \rangle \equiv C \langle x_0,\ldots,x_{n-2},0,0 \rangle $$
and $E$ is standardly semilinear. We will use the following notation for retracts of $C$
(the table illustrates their mutual arrangement,
where the last and last-but-one coordinates of $\Sigma^{n+1}$
are thought as ordinate and abscissa respectively; for example, $B_2$ corresponds
to fixing abscissa by $2$):

\noindent\mbox{}\hfill$
\begin{array}{rcl}
  E_{a,b}\langle x_0,\ldots,x_{n-2} \rangle & \triangleq & C \langle x_0,\ldots,x_{n-2},a,b \rangle, \\
  A_b\langle x_0,\ldots,x_{n-2},y \rangle & \triangleq & C\langle x_0,\ldots,x_{n-2},y,b \rangle, \\
  B_a\langle x_0,\ldots,x_{n-2},z \rangle & \triangleq & C\langle x_0,\ldots,x_{n-2},a,z \rangle.
\end{array}
$
\raisebox{-17.5mm}{\small
\def\dick{1pt}
\def\normal{0.4pt}
\def\duenn{0.1pt}
\unitlength 0.7mm
\def\emline#1#2#3#4#5#6{%
\put(#1,#2){\special{em:moveto}}%
\put(#4,#5){\special{em:lineto}}}
\begin{picture}(52.00,49.00)
\special{em:linewidth \normal}
\emline{ 9.00}{47.00}{1}{17.00}{47.00}{2}
\emline{17.00}{39.00}{1}{ 9.00}{39.00}{2}
\emline{ 9.00}{39.00}{1}{ 9.00}{47.00}{2}
\emline{19.00}{47.00}{1}{27.00}{47.00}{2}
\emline{27.00}{47.00}{1}{27.00}{39.00}{2}
\emline{29.00}{47.00}{1}{37.00}{47.00}{2}
\emline{37.00}{47.00}{1}{37.00}{39.00}{2}
\emline{37.00}{39.00}{1}{29.00}{39.00}{2}
\emline{29.00}{39.00}{1}{29.00}{47.00}{2}
\emline{39.00}{47.00}{1}{47.00}{47.00}{2}
\emline{47.00}{39.00}{1}{39.00}{39.00}{2}
\emline{39.00}{39.00}{1}{39.00}{47.00}{2}
\emline{17.00}{37.00}{1}{17.00}{29.00}{2}
\emline{17.00}{29.00}{1}{ 9.00}{29.00}{2}
\emline{ 9.00}{29.00}{1}{ 9.00}{37.00}{2}
\emline{19.00}{37.00}{1}{27.00}{37.00}{2}
\emline{27.00}{37.00}{1}{27.00}{29.00}{2}
\emline{27.00}{29.00}{1}{19.00}{29.00}{2}
\emline{19.00}{29.00}{1}{19.00}{37.00}{2}
\emline{29.00}{37.00}{1}{37.00}{37.00}{2}
\emline{37.00}{37.00}{1}{37.00}{29.00}{2}
\emline{37.00}{29.00}{1}{29.00}{29.00}{2}
\emline{29.00}{29.00}{1}{29.00}{37.00}{2}
\emline{39.00}{37.00}{1}{47.00}{37.00}{2}
\emline{47.00}{37.00}{1}{47.00}{29.00}{2}
\emline{47.00}{29.00}{1}{39.00}{29.00}{2}
\emline{39.00}{29.00}{1}{39.00}{37.00}{2}
\emline{ 9.00}{27.00}{1}{17.00}{27.00}{2}
\emline{17.00}{27.00}{1}{17.00}{19.00}{2}
\emline{17.00}{19.00}{1}{ 9.00}{19.00}{2}
\emline{ 9.00}{19.00}{1}{ 9.00}{27.00}{2}
\emline{19.00}{27.00}{1}{27.00}{27.00}{2}
\emline{27.00}{27.00}{1}{27.00}{19.00}{2}
\emline{27.00}{19.00}{1}{19.00}{19.00}{2}
\emline{19.00}{19.00}{1}{19.00}{27.00}{2}
\emline{29.00}{27.00}{1}{37.00}{27.00}{2}
\emline{37.00}{27.00}{1}{37.00}{19.00}{2}
\emline{37.00}{19.00}{1}{29.00}{19.00}{2}
\emline{29.00}{19.00}{1}{29.00}{27.00}{2}
\emline{39.00}{27.00}{1}{47.00}{27.00}{2}
\emline{47.00}{27.00}{1}{47.00}{19.00}{2}
\emline{47.00}{19.00}{1}{39.00}{19.00}{2}
\emline{39.00}{19.00}{1}{39.00}{27.00}{2}
\emline{ 9.00}{17.00}{1}{17.00}{17.00}{2}
\emline{17.00}{17.00}{1}{17.00}{ 9.00}{2}
\emline{17.00}{ 9.00}{1}{ 9.00}{ 9.00}{2}
\emline{ 9.00}{ 9.00}{1}{ 9.00}{17.00}{2}
\emline{19.00}{17.00}{1}{27.00}{17.00}{2}
\emline{27.00}{17.00}{1}{27.00}{ 9.00}{2}
\emline{27.00}{ 9.00}{1}{19.00}{ 9.00}{2}
\emline{19.00}{ 9.00}{1}{19.00}{17.00}{2}
\emline{29.00}{17.00}{1}{37.00}{17.00}{2}
\emline{37.00}{17.00}{1}{37.00}{ 9.00}{2}
\emline{37.00}{ 9.00}{1}{29.00}{ 9.00}{2}
\emline{29.00}{ 9.00}{1}{29.00}{17.00}{2}
\emline{39.00}{17.00}{1}{47.00}{17.00}{2}
\emline{47.00}{17.00}{1}{47.00}{ 9.00}{2}
\emline{47.00}{ 9.00}{1}{39.00}{ 9.00}{2}
\emline{39.00}{ 9.00}{1}{39.00}{17.00}{2}
\emline{19.00}{47.00}{1}{19.00}{39.00}{2}
\emline{19.00}{39.00}{1}{27.00}{39.00}{2}
\emline{17.00}{47.00}{1}{17.00}{39.00}{2}
\emline{17.00}{37.00}{1}{ 9.00}{37.00}{2}
\emline{47.00}{47.00}{1}{47.00}{39.00}{2}
\emline{ 8.00}{49.00}{1}{48.00}{49.00}{2}
\emline{48.00}{49.00}{1}{49.00}{48.00}{2}
\emline{49.00}{48.00}{1}{49.00}{ 8.00}{2}
\emline{49.00}{ 8.00}{1}{48.00}{ 7.00}{2}
\emline{48.00}{ 7.00}{1}{ 8.00}{ 7.00}{2}
\emline{ 8.00}{ 7.00}{1}{ 7.00}{ 8.00}{2}
\emline{ 7.00}{ 8.00}{1}{ 7.00}{48.00}{2}
\emline{ 7.00}{48.00}{1}{ 8.00}{49.00}{2}
\emline{28.00}{47.00}{1}{29.00}{48.00}{2}
\emline{29.00}{48.00}{1}{37.00}{48.00}{2}
\emline{37.00}{48.00}{1}{38.00}{47.00}{2}
\emline{38.00}{47.00}{1}{38.00}{ 6.00}{2}
\emline{38.00}{ 6.00}{1}{37.00}{ 5.00}{2}
\emline{37.00}{ 5.00}{1}{29.00}{ 5.00}{2}
\emline{29.00}{ 5.00}{1}{28.00}{ 6.00}{2}
\emline{28.00}{ 6.00}{1}{28.00}{47.00}{2}
\emline{ 9.00}{18.00}{1}{ 8.00}{19.00}{2}
\emline{ 8.00}{19.00}{1}{ 8.00}{27.00}{2}
\emline{ 8.00}{27.00}{1}{ 9.00}{28.00}{2}
\emline{ 9.00}{28.00}{1}{50.00}{28.00}{2}
\emline{50.00}{28.00}{1}{51.00}{27.00}{2}
\emline{51.00}{27.00}{1}{51.00}{19.00}{2}
\emline{51.00}{19.00}{1}{50.00}{18.00}{2}
\emline{50.00}{18.00}{1}{ 9.00}{18.00}{2}
\emline{-1.00}{ 5.00}{1}{-1.00}{13.00}{2}
\emline{-1.00}{13.00}{1}{-0.51}{11.06}{2}
\emline{-1.49}{11.06}{1}{-1.00}{13.00}{2}
\emline{-1.00}{ 5.00}{1}{ 7.00}{ 5.00}{2}
\emline{ 7.00}{ 5.00}{1}{ 5.06}{ 4.51}{2}
\emline{ 5.06}{ 5.49}{1}{ 7.00}{ 5.00}{2}
\put(13.00,43.00){\makebox(0,0)[cc]{$E_{0,3}$}}
\put(23.00,43.00){\makebox(0,0)[cc]{$E_{1,3}$}}
\put(33.00,43.00){\makebox(0,0)[cc]{$E_{2,3}$}}
\put(43.00,43.00){\makebox(0,0)[cc]{$E_{3,3}$}}
\put(13.00,33.00){\makebox(0,0)[cc]{$E_{0,2}$}}
\put(23.00,33.00){\makebox(0,0)[cc]{$E_{1,2}$}}
\put(33.00,33.00){\makebox(0,0)[cc]{$E_{2,2}$}}
\put(43.00,33.00){\makebox(0,0)[cc]{$E_{3,2}$}}
\put(13.00,23.00){\makebox(0,0)[cc]{$E_{0,1}$}}
\put(23.00,23.00){\makebox(0,0)[cc]{$E_{1,1}$}}
\put(33.00,23.00){\makebox(0,0)[cc]{$E_{2,1}$}}
\put(43.00,23.00){\makebox(0,0)[cc]{$E_{3,1}$}}
\put(13.00,13.00){\makebox(0,0)[cc]{$E_{0,0}$}}
\put(23.00,13.00){\makebox(0,0)[cc]{$E_{1,0}$}}
\put(33.00,13.00){\makebox(0,0)[cc]{$E_{2,0}$}}
\put(43.00,13.00){\makebox(0,0)[cc]{$E_{3,0}$}}
\put( 7.00, 4.00){\makebox(0,0)[tr]{$y$}}
\put(-2.00,13.00){\makebox(0,0)[tr]{$z$}}
\put(33.00, 4.00){\makebox(0,0)[tc]{$B_2$}}
\put(52.00,23.00){\makebox(0,0)[cl]{$A_1$}}
\put(49.00, 7.00){\makebox(0,0)[tl]{$C$}}
\put(52.00,13.00){\makebox(0,0)[cl]{$A_0$}}
\put(52.00,33.00){\makebox(0,0)[cl]{$A_2$}}
\put(52.00,43.00){\makebox(0,0)[cl]{$A_3$}}
\put(13.00, 4.00){\makebox(0,0)[tc]{$B_0$}}
\put(23.00, 4.00){\makebox(0,0)[tc]{$B_1$}}
\put(43.00, 4.00){\makebox(0,0)[tc]{$B_3$}}
\end{picture}
}
\hfill\mbox{}

Since $A_0$ is reducible and fixing $y:=0$ leads to the irreducible $E$,
by Proposition~\ref{p:canon} we have
\begin{equation}\label{eq:A0}
A_0\langle x_0,\ldots,x_{n-2},y \rangle
\equiv
E\langle x_0,\ldots,x_{i-1},h(x_i,y),x_{i+1},\ldots,x_{n-2} \rangle
\end{equation}
for some $i\in \{0,\ldots n-2\} $ and $2$-qua\-si\-group $h$ such that $h(x_i,0)\equiv x_i$.

From (\ref{eq:A0}), we see that all the retracts
$E_{a,0}$, $a\in\Sigma$, are isotopic to $E$. Similarly, we can get the following:

\emph{Claim~(*)}. All the retracts $E_{a,b}$, $a,b\in\Sigma$ are isotopic to $E$.

Then,
we conclude that a representation similar to (\ref{eq:A0}) is valid for every $b\in\Sigma$:
$$A_b \langle x_0,\ldots,x_{n-2},y \rangle
\equiv
E_{0,b} \langle x_0,\ldots,x_{i_b-1},h_b(x_{i_b},y),x_{i_b+1},\ldots,x_{n-2} \rangle
$$
for some $i\in \{0,\ldots n-2\} $ and $2$-qua\-si\-group $h_b$ such that $h_b(x,0)\equiv x$.

\emph{Claim~(**)}. We claim that $i_b$ does not depend on $b$.
Indeed, assume, for example, that $i_1=0$ and $i_2=1$, i.\,e.,
\begin{eqnarray*}
A_1 \langle x_0,\ldots,x_{n-2},y \rangle
&\equiv&
E_{0,1} \langle h_1(x_0,y),x_1,x_2,\ldots,x_{n-2} \rangle,\\
A_2 \langle x_0,\ldots,x_{n-2},y \rangle
&\equiv&
E_{0,2} \langle x_0,h_2(x_1,y),x_2,\ldots,x_{n-2} \rangle.
\end{eqnarray*}
Then, fixing $x_0$ in the first case leads to a retract isotopic to $E$;
fixing $x_0$ in the second case leads to a reducible retract (recall that $n\geq 5$).
But, analogously to (*), these two retracts are isotopic;
this contradicts the irreducibility of $E$ and proves (**).

Without loss of generality we can assume that $i_b=0$, i.\,e.,
\begin{equation}\label{eq:A0E}
A_b \langle x_0,x_1,\tilde x_2,y \rangle \equiv E_{0,b} \langle
h_b(x_0,y),x_1,\tilde x_2 \rangle;
\end{equation}
here and later $\tilde x_2\triangleq (x_2,\ldots,x_{n-2})$. Similarly, we can assume without loss of generality that either
\begin{equation}\label{eq:B0E}
B_a \langle x_0,x_1,\tilde x_2,z \rangle \equiv E_{a,0} \langle
g_a(x_0,z),x_1,\tilde x_2 \rangle
\end{equation}
or
\begin{equation}\label{eq:B1E}
B_a \langle x_0,x_1,\tilde x_2,z \rangle \equiv E_{a,0} \langle
x_0,g_a(x_1,z),\tilde x_2 \rangle
\end{equation}
where $2$-qua\-si\-groups $g_a$ satisfy $g_a(x,0)\equiv x$.

Using (\ref{eq:A0E}) and (\ref{eq:B0E}), we derive
\begin{eqnarray}
C \langle x_0,x_1,\tilde x_2,y,z\rangle
&\equiv& \nonumber
A_z\langle x_0,x_1,\tilde x_2,y \rangle
\\&\equiv& \nonumber
E_{0,z}\langle h_z(x_0,y),x_1,\tilde x_2\rangle
\\&\equiv& \label{eq:4-5}
B_{0}\langle h_z(x_0,y),x_1,\tilde x_2, z\rangle
\\&\equiv& \nonumber
E_{0,0} \langle g_0(h_z(x_0,y),z),x_1,\tilde x_2 \rangle,
\end{eqnarray}
which means that $C$ is reducible, because
$f(x,y,z) \triangleq g_0(h_z(x,y),z)$ must be a $3$-qua\-si\-group.
So, it remains to consider the case (\ref{eq:B1E}).
Consider two subcases.

\emph{Case}~1. The $2$-qua\-si\-group $g_a$ does not depend on $a$;
denote $g \triangleq g_a$.
Then, repeating the first three steps of (\ref{eq:4-5})
and applying (\ref{eq:B1E}), we derive that
$$
C\langle x_0,x_1,\tilde x_2,y,z \rangle \equiv
E_{0,0}\langle h_0(x_0,y), g(x_1,z), \tilde x_2 \rangle,
$$
and $C$ is reducible.

\emph{Case}~2. For some fixed $a$ we have $g_0\neq g_a$;
denote $s_i(x)\triangleq g_0(x,i)$, $t_i(x)\triangleq g_a(x,i)$,
and $r_i(x) \triangleq h_i(x,a)$.
From (\ref{eq:A0E}), we see that
\begin{eqnarray}
E_{a,0}\langle x_0,x_1,\tilde x_2 \rangle
&\equiv & \label{eq:E-1}
E_{0,0}\langle r_0(x_0),x_1,\tilde x_2 \rangle,
\\
E_{a,b}\langle x_0,x_1,\tilde x_2 \rangle
&\equiv & \label{eq:E-2}
E_{0,b}\langle r_b(x_0),x_1,\tilde x_2 \rangle.
\end{eqnarray}
From (\ref{eq:B1E}), we see that
\begin{eqnarray}
E_{0,b}\langle x_0,x_1,\tilde x_2 \rangle
&\equiv & \label{eq:E-3}
E_{0,0}\langle x_0,s_b(x_1),\tilde x_2 \rangle,
\\
E_{a,b}\langle x_0,x_1,\tilde x_2 \rangle
&\equiv & \label{eq:E-4}
E_{a,0}\langle x_0,t_b(x_1),\tilde x_2 \rangle.
\end{eqnarray}
Applying consequtively (\ref{eq:E-3}), (\ref{eq:E-2}), (\ref{eq:E-4}), and (\ref{eq:E-1}),
we find that for each $b$ the retract  $E=E_{0,0}$ satisfies
$$
E\langle x_0,x_1,\tilde x_2 \rangle
\equiv E_{0,b}\langle\ldots\rangle
\equiv E_{a,b}\langle\ldots\rangle
\equiv E_{a,0}\langle\ldots\rangle
\equiv E\langle r_0 r_b^{-1} x_0,t_b s_b^{-1} x_1,\tilde x_2 \rangle.
$$
By Proposition~\ref{p:aut}, the irreducibility of $E$ means that
$t_b s_b^{-1} \in\{Id,(1{,}0{,}3{,}2)\}$
for every $b$.
By Proposition~\ref{p:2-2}, for some permutation $\phi$
the $2$-qua\-si\-group $s(x,z)\triangleq g_0( x,\phi z)$ is standardly semilinear.
Since a composition of standardly semilinear qua\-si\-groups is a standardly semilinear qua\-si\-group,
we see that $B_0$ is a semilinear $(n-1)$-qua\-si\-group. Lemma~\ref{pro:3b} completes the proof.%
\qquad\endproof

\section{Acknowledgement}
The authors thank the referees for their work in reviewing the manuscript
and the audience of the seminar ``Coding Theory''
in the Sobolev Institute of Mathematics
for the patience during reporting this result.



\providecommand\href[2]{#2} \providecommand\url[1]{\href{#1}{#1}}

\def\DOI#1{{\tt{DOI}: \href{http://dx.doi.org/#1}{#1}}}
\def\ArXiv#1{{\tt{ArXiv}: \href{http://arXiv.org/abs/#1}{#1}}}

\end{document}